\documentstyle[amscd,11pt]{amsart}
\author{Kevin P. Knudson}\thanks{Supported by an NSF Postdoctoral Fellowship, grant no. DMS-9627503}
\address{Department of Mathematics, Northwestern University, 
Evanston, IL  60208}
\email{knudson@@math.nwu.edu}
\date{December 4, 1997}
\title[Linear Groups over Hensel Rings]{Low Dimensional Homology
of Linear Groups over Hensel Local Rings}

\newcommand{\slnk}{SL_n(k)}

\newcommand{\sln}{SL_n(k[[t_1,t_2,\dots ,t_m]])}
\newcommand{\glnk}{GL_n(k)}
\newcommand{\gln}{GL_n(R)}

\newcommand{\pglk}{PGL_2(k)}
\newcommand{\pgl}{PGL_2(R)}
\newcommand{\zz}{{\Bbb Z}}
\newcommand{\zp}{{\Bbb Z}/p}

\newcommand{\lra}{\longrightarrow}

\newcommand{\sps}{spectral sequence }
\newcommand{\uk}{k^{\times}}

\newcommand{\ur}{R^{\times}}

\newcommand{\cnm}{C_{n,m}}
\newcommand{\lnm}{L_{n,m}}
\newcommand{\bop}{\bigoplus}
\newcommand{\ra}{\rightarrow}

\newtheorem{theorem}{Theorem}[section]
\newtheorem{prop}[theorem]{Proposition}
\newtheorem{lemma}[theorem]{Lemma}
\newtheorem{cor}[theorem]{Corollary}
\newtheorem{conj}[theorem]{Conjecture}

\theoremstyle{remark}

\begin{document}

\begin{abstract}  We prove that if $R$ is a Hensel local ring with
infinite residue field $k$, the natural map $H_i(GL_n(R),\zp) \ra H_i(GL_n(k),\zp)$ is an isomorphism for $i\le 3$, $p \ne \text{char}\; k$.
This implies rigidity for $H_i(GL_n)$, $i \le 3$, which in turn implies the Friedlander--Milnor conjecture in positive characteristic in degrees
$\le 3$.
\end{abstract}

\maketitle

A fundamental question in the homology of linear groups is that of
rigidity:  given a smooth affine curve $X$ over an algebraically closed
field $k$ and closed points $x,y$ on $X$, do the corresponding 
specialization homomorphisms
$$s_x,s_y:H_*(G(k[X]),\zp) \lra H_*(G(k),\zp)$$
coincide?  Here, $G$ is a reductive algebraic group and $p$ is a prime
not equal to the characteristic of $k$.  The answer is yes when $X$ is
the affine line and $G=SL_n, GL_n, PGL_n$ since the inclusion $G(k)\ra G(k[t])$
induces an isomorphism
$$H_*(G(k),\zz) \lra H_*(G(k[t]),\zz)$$
(see \cite{knudson}) and the map is split by evaluation at any $x \in {\Bbb A}^1$.

Rigidity in algebraic $K$-theory has spectacular consequences, including
the calculation of the $K$-theory of algebraically closed fields and the
solution of the Friedlander--Milnor conjecture for the stable general
linear group $GL$.  Similarly, a proof of rigidity for $G(k[X])$ would
imply the Friedlander--Milnor conjecture for $G$ (see Section \ref{fmconj}).

Rigidity would follow if one could prove the following stronger result.
Let $X$ be a smooth curve over an algebraically closed field $k$ and let
$x$ be a closed point on $X$.  Denote by ${\cal O}_x^h$ the henselization of
${\cal O}_x$.

\medskip

\noindent {\bf Conjecture.} {\em The inclusion $G(k) \ra G({\cal O}_x^h)$
induces an isomorphism
$$H_*(G(k),\zp) \lra H_*(G({\cal O}_x^h),\zp)$$
for $p \ne \text{\em char}\; k$.}

\medskip

In Section \ref{fmconj} we sketch a proof of rigidity assuming this conjecture (this argument is well-known).  Our results will allow us to
deduce the Friedlander--Milnor conjecture in positive characteristic
for $H_i(GL_n)$, $i\le 3$.  The conjecture is
known to be true for
$H_3(SL_2({\Bbb C}))$ by the work of Sah \cite{sah}.

In this note we study $H_*(G(R),\zp)$ for $G=SL_n$
and $G=GL_n$, where $R$ is a Hensel local $k$-algebra and
$k$ is an infinite field.  Our results are far from complete.  We prove the following
result.

\medskip

\noindent {\bf Theorem.}  {\em The inclusion $\glnk \ra \gln$ induces an
isomorphism
$$H_i(\glnk,\zp) \lra H_i(\gln,\zp)$$
for $i\le 3$.}

\medskip

Note that this map is split injective for all $i$; the only issue is
surjectivity.

The expert reader by now will have noticed that this theorem contains
only one new case, namely $H_3(GL_2)$.  The above map is in fact an
isomorphism for $i\le n$ by a combination of Suslin's stability theorem
\cite{nessus} and rigidity in $K$-theory.  However, we show a bit more
than the theorem states.   We first prove that the map
$$H_2(\pglk,\zz/p) \lra H_2(\pgl,\zz/p)$$
is an isomorphism and use this to deduce the result for $H_2(GL_2)$.
Moreover, our approach is different from that used by Suslin \cite{suslin1} to compute $H_\bullet(GL(R),\zz/p)$.  We then study Bloch
groups to obtain the result for $H_3(GL_2)$.

The corresponding statement for $G=SL_n$ and $i=0,1$ is almost obvious.  However, we provide an alternate
proof for $i=1$ in a special case as follows.  Let $F$ be any field and denote by ${\frak m}$ the maximal ideal $(t_1,t_2,\dots ,t_m)$ of 
$F[[t_1,t_2,\dots ,t_m]]$.  Consider the short exact sequence

\begin{equation}\label{ss1}
1 \lra \cnm \lra SL_n(F[[t_1,t_2,\dots ,t_m]]) \lra SL_n(F) \lra 1.
\end{equation}
Here, 
the group $\cnm$ consists of those matrices which are congruent to
the identity modulo ${\frak m}$.  We have the following result.

\medskip

\noindent {\bf Proposition.}  {\em If $n\ge 3$, then $$H_1(\cnm,\zz) = 
\displaystyle \bop_{i=1}^m {\frak sl}_n(F).$$  When $n=2$, assume
further that $\text{\em char}\; F \ne 2$ and $F\ne {\Bbb F}_3$.  Then $$H_1(C_{2,m},\zz)
= \displaystyle \bop_{i=1}^m {\frak sl}_2(F).$$}

\medskip

When $p$ is a prime distinct from the characteristic
of $k$, the proposition implies that $H_1(\cnm,\zp) = 0$ unless
$n=2, \text{char}\; k = 2$ or $F={\Bbb F}_3$.  In this case, however, one can still conclude
that $H_1(C_{2,m},\zp)=0$ by other methods.  One expects that
$H_\bullet(C_{n,m},\zz/p)$ vanishes in all positive degrees.  We provide evidence for this in Section 2. 

\medskip

\noindent {\em Conventions.}  Throughout, $k$ denotes a field and $p$ is
a prime distinct from the characteristic of $k$.

\section{The Congruence Subgroup}\label{congruence}

Recall that the group
$\cnm$ is the kernel of the natural map
$$SL_n(F[[t_1,t_2,\dots ,t_m]]) \stackrel{\text{mod}\; {\frak m}}{\lra}
SL_n(F)$$
where ${\frak m}$ is the ideal $(t_1,t_2,\dots ,t_m)$.  Define a sequence
of subgroups $\cnm^i$ by
$$\cnm^i=\{X \in SL_n(F[[t_1,t_2,\dots ,t_m]]): X \equiv I\;\text{mod}\; 
{\frak m}^{i} \}.$$
If no confusion can result, we usually drop the subscripts from the 
notation.  Note that each $C^i$ is a normal subgroup as it is the kernel of the map
$$SL_n(F[[t_1,t_2,\dots ,t_m]]) \lra SL_n(F[[t_1,t_2,\dots ,t_m]]/ 
{\frak m}^i).$$

For each $i$, define a homomorphism
$$\rho_i: \cnm^i \lra \bop_{\lambda} {\frak sl}_n(F)$$
by
$$\rho_i(I+\sum_{\lambda} t_1^{k_1}t_2^{k_2}\cdots
t_m^{k_m} X_{\lambda} + \cdots ) = (X_{\lambda})_{\lambda}$$
where $\lambda=(k_1,k_2,\dots ,k_m)$ is an $m$-partition of $i$ ({\em
i.e.,} $\lambda$ is an $m$-tuple of nonnegative integers whose sum is
$i$)
and $(X_{\lambda})_{\lambda}$ denotes the 
element of $\bop {\frak sl}_n(F)$ given by the various $X_{\lambda}$.  Note
that $\rho_i$ is well-defined since the equation
$$1=\det X \equiv 1 + \sum_{\lambda} t_1^{k_1}\cdots
t_m^{k_m} \text{tr} X_{\lambda}\; \text{mod}\; {\frak m}^{i+1}$$
implies that $\text{tr}\, X_{\lambda} = 0$ for each partition $\lambda$.  One checks
easily that $\rho_i$ is a surjective group homomorphism with kernel
$\cnm^{i+1}$ (see \cite{lee} for an analogous result for congruence 
subgroups of $SL_n(\zz)$).  Hence, for each $i$ we have an isomorphism
$$\cnm^i/\cnm^{i+1} \cong \bop_{\lambda} {\frak sl}_n(F).$$

\medskip

\noindent {\em Remark.}  The above discussion remains valid for any
simple algebraic group $G$.  The quotients $C^i/C^{i+1}$ are direct sums of
copies of the lie algebra ${\frak g}$.  We have chosen to work with $SL_n$ just to
fix ideas.

\medskip

\begin{lemma}\label{l1}  For each $i,j$, $[C^i,C^j] \subseteq C^{i+j}$.
\end{lemma}

\begin{pf}  This follows easily by writing out $XYX^{-1}Y^{-1}$ for
$X\in C^i$, $Y\in C^j$ and noting that the inverse of $I+\sum Z_{\lambda} +\cdots$
is $I-\sum Z_{\lambda}+\cdots$.
\end{pf}

For any group $G$, denote by $\Gamma^\bullet$ the lower central series
of $G$.

\begin{cor} For each $i$, $\Gamma^i \subseteq \cnm^i$.
\end{cor}

\begin{pf}  The series $\cnm^\bullet$ is a descending central series and
as such contains the lower central series.
\end{pf}

We now recall the following theorem of Klingenberg \cite{kling}.

\begin{theorem}  If $A$ is a local ring then for $n\ge 3$ the only normal
subgroups of $SL_n(A)$ are the congruence subgroups.  If $n=2$, the same is true as long as the residue characteristic of $A$ is not $2$ or the residue field is not ${\Bbb F}_3$. \hfill $\qed$
\end{theorem}

\begin{cor}  The filtration $\cnm^\bullet$ is the lower 
central series. 
\end{cor}

\begin{pf}  Note that the groups $\Gamma^i$ are characteristic subgroups
of $\cnm$ and hence are normal in $SL_n(F[[t_1,\dots ,t_m]])$.  Since there are clearly elements in $\Gamma^i - \cnm^{i+1}$, it follows that
$\Gamma^i = \cnm^i$ for all $i$.
\end{pf}

\begin{cor}  For all $n\ge 3$, $H_1(\cnm,\zz) = \bop_{i=1}^m {\frak sl}_n(F)$.  If $\text{\em char}\; F \ne 2$ or $F\ne {\Bbb F}_3$, then $H_1(C_{2,m},\zz)
= \bop_{i=1}^m {\frak sl}_2(F)$.
\end{cor}

\noindent {\em Proof.}  Since $\cnm^2 = \Gamma_{n,m}^2$, we have
$$H_1(\cnm,\zz) = \cnm/\Gamma^2=\cnm/\cnm^2 = \bop_{i=1}^m
{\frak sl}_n(F).$$ \hfill $\qed$

\noindent {\em Remark.}  This corollary holds for power series rings over more general rings (for example, rings of the form $\zz[J^{-1}]$ for a set of primes $J$).

\begin{cor}  If $p$ is a prime different from the
characteristic of $F$, then
$$H_1(\cnm,\zp) = 0.$$
\end{cor}

\begin{pf}  This follows from the previous corollary except in the case
$n=2$, $\text{char}\; F = 2$ or $F={\Bbb F}_3$.  However, one could also prove this corollary
directly by noting that each element of $\cnm$ has a $p$th root in
$\cnm$ and hence that $H_1(\cnm,\zp) = H_1(\cnm,\zz)\otimes
\zp = 0$.
\end{pf}

\begin{cor}  For each $n$,
$$H_1(SL_n(F[[t_1,\dots ,t_m]]),\zp) \cong H_1(SL_n(F),\zp).$$
\end{cor}

\begin{pf}  This follows easily by considering the
Hochschild--Serre spectral sequence associated to the extension
(\ref{ss1}).
\end{pf}

\section{Conjectures about $H_*(\cnm,\zp)$}\label{conjectures}

The isomorphism $$H_\bullet(\sln,\zp) \cong H_\bullet(\slnk,\zp)$$ would follow
easily if one could prove the following.

\begin{conj}  For all $l>0$, $H_l(\cnm(k[[t_1,t_2,\dots ,t_m]]),\zp) = 0$.
\end{conj}

We have just proved the case $l=1$.  We make the following observations.

\subsection{Continuous homology}  For each $i\ge 1$, denote by $\lnm^i$ the kernel of the homomorphism
$$SL_n(k[t_1,\dots ,t_m]) \lra SL_n(k[t_1,\dots ,t_m]/(t_1,\dots t_m)^i).$$  The group $\cnm$ is the inverse
limit of the nilpotent groups $\lnm/\lnm^i$.  One checks that these
groups have no homology with $\zp$ coefficients by noting that the
successive graded quotients $\lnm^j/\lnm^{j+1}$ of
$\lnm/\lnm^i$ are $k$-vector spaces (see for example \cite{suslin1});
an iterated use of the Hochschild--Serre \sps then shows that $\lnm/\lnm^i$
is $\zp$--acyclic.  It follows that the ``continuous homology"
$$\varprojlim H_\bullet(\cnm/\cnm^i,\zp) = \varprojlim H_\bullet(\lnm/\lnm^i,\zp)$$
vanishes and that
$$H_\bullet(SL_n(k[t_1,t_2,\dots ,t_m]/{\frak m}^i),\zp) \cong H_\bullet(\slnk,\zp)$$ for all $i\ge 1$.

\subsection{Large acyclic subgroups}  When $m=1$, it is easy to see
that $$H_\bullet(L_{n,1},\zp) = 0$$ (see \cite{knudson}) and the group $L_{n,1}$
is a subgroup of $C_{n,1}$.

\subsection{Product spaces}  The group $\cnm$ is the subgroup of
$\prod_i \lnm/\lnm^i$ consisting of coherent sequences.  A theorem of
P. Goerss \cite{goerss} implies that this product is $\zp$--acyclic.

\subsection{Cosimplicial replacement}  We have seen that $\cnm$ is the
inverse limit of the nilpotent groups $\lnm/\lnm^i$.  Let $X_i$ be a
$K(\lnm/\lnm^i,1)$ space and consider the tower of fibrations
$$X_2 \leftarrow X_3 \leftarrow X_4 \leftarrow \cdots$$
The homotopy inverse limit of this tower is a $K(\cnm,1)$-space;
denote it by $X$.  Computing the homology of homotopy inverse limits
is an important problem in homotopy theory, and has been studied 
intensely by Bousfield \cite{bousfield}, Goerss \cite{goerss},
Shipley \cite{shipley}, and others.  

The standard technique is to consider the {\em cosimplicial replacement}
of the tower.  This is a cosimplicial space $X^\bullet$ whose $n$th space
$X^n$ is a certain product of the spaces $X_i$ of the tower (see Bousfield--Kan \cite[XI, Sec. 5]{bk}).  Given a cosimplicial space $X^\bullet$, there
is an associated homology spectral sequence \cite{bousfield}
$$E_{m,t}^2 = \pi^mH_t(X^\bullet;A)$$
for any abelian group $A$.  Here, the term on the right is the $m$th
cohomotopy of the cosimplicial abelian group $H_t(X^\bullet;A)$.  The
main problem is that of convergence of this spectral sequence.  Under
certain conditions, the \sps converges to the homology of the total space
of $X^\bullet$.  In the case of the cosimplicial replacement of a tower
of fibrations, the total space is the homotopy inverse limit of the tower.

In the case of $\cnm$, we see that the cosimplicial replacement is
particularly nice.  Its spaces $X^n$ are certain products of the $X_i$ and
by Goerss' theorem each of these is $\zp$--acyclic.  It follows that the
associated homology spectral sequence has $E^2$--term
$$E_{m,t}^2 = \begin{cases}
                            \zp      &   (m,t)=(0,0)   \\
                              0       &   (m,t)\ne (0,0).
                     \end{cases}$$
The question is what, if anything, is this spectral sequence converging to?
Unfortunately, all known convergence theorems do not apply to the situation at hand.  Most of these require some hypotheses on the spaces
$X^n$ such as nilpotency (which we do not have in this case---we would
need the groups $\lnm/\lnm^i$ to have bounded nilpotence degree).  A
proof that the \sps converges to the homology of the total space of $X^\bullet$ in this 
case would immediately show that $\cnm$ is $\zp$--acyclic.

\medskip

Given all of this, it is difficult to imagine that $\cnm$ has any $\zp$
homology.  Speaking heuristically, the group $\cnm$ is built up by 
successively attaching copies of ${\frak sl}_n(k)$ and since this
group is $\zp$--acyclic, one would expect that $\cnm$ is as well.

\section{Spectral Sequences and the Second Homology Groups}\label{h2}
Hereafter, $R$ denotes a Hensel local $k$-algebra, and 
${\frak m}$ denotes the maximal ideal of $R$.  We have $R/{\frak m}
= k$.

We begin by considering the group $\pgl$.  We will show that the
inclusion $\pglk \ra \pgl$ induces an isomorphism on second homology
with coefficients in $\zp$.

Let $A$ be a local ring with maximal ideal ${\frak m}$ and residue field
$F$, which we assume to be infinite.  Recall that a column vector $v\in
A^2$ is {\em unimodular} if the ideal generated by its entries is $A$
itself.  Denote by $\overline{v}$ the vector in $F^2$ obtained by reducing
the entries of $v$ modulo ${\frak m}$.  We say that a collection 
$v_1,v_2,\dots ,v_k$ of unimodular vectors is in {\em general position}
if the collection $\overline{v}_1,\overline{v}_2,\dots ,\overline{v}_k$
is in general position ({\em i.e.,} the matrix determined by any pair of
them is invertible).

Construct a simplicial set $S_\bullet$ as follows.  The nondegenerate
$p$--simplices are collections $v_0,v_1,\dots ,v_p$ of projective
equivalence classes of unimodular vectors which are in general position
(this amounts to saying that the $v_i$ are distinct closed
points in ${\Bbb P}^1$).  Denote by $C_\bullet(A)$ the associated 
simplicial chain complex.  A proof of the following may be found in
Nesterenko--Suslin \cite{nessus}.

\begin{lemma}  The augmented complex $C_\bullet(A) \stackrel{\epsilon}{\ra}
\zz \ra 0$ is acyclic.\hfill $\qed$
\end{lemma}

Note that $PGL_2(A)$ acts transitively on $C_i(A)$ for $i\le 2$ ({\em i.e.,}
$PGL_2(A)$ acts 3-transitively on points in ${\Bbb P}^1(A)$).  Denote by
$0$ the vector $\left(\begin{array}{c}
                                          1   \\   0
                                   \end{array}\right)$;
by $\infty$ the vector $\left(\begin{array}{c}
                                          0   \\   1
                                   \end{array}\right)$;
and by $1$ the vector $\left(\begin{array}{c}
                                          1   \\   1
                                   \end{array}\right)$.
Then each orbit of the action of $PGL_2(A)$ on $C_p(A)$ has a unique
representative of the form $(0,\infty,1,v_1,\dots ,v_{p-2})$.  Since
each $v_i$ is in general position with $0,\infty$ and $1$, we see that both
entries of $v_i$ must be units in $A$.  Since we are using projective
equivalence classes of vectors, we may then assume that the vector
$v_i$ has the form $\left(\begin{array}{c}
                                          1   \\   \alpha_i
                                   \end{array}\right)$
where $\alpha_i$ is a unit in $A$.  Moreover, the $\alpha_i$ satisfy the
additional conditions that $1-\alpha_i \in A^\times$ (since $v_i$ is
in general position with $1$) and $\alpha_i - \alpha_j \in A^\times$
for $i\ne j$ (since $v_i$ is in general position with $v_j$).

Denote by $E_{\bullet,\bullet}^\bullet(A)$ the spectral sequence 
associated to the action of the group $PGL_2(A)$ on $C_\bullet(A)$.  This \sps converges to the homology of $PGL_2(A)$
and has $E^1$--term
$$E_{p,q}^1 = \bop_{\sigma \in \Sigma_p} H_q(G_\sigma)$$
where $\Sigma_p$ is a set of representatives of the orbits of the action
of $PGL_2(A)$ on $C_p(A)$ and $G_\sigma$ is the stabilizer of $\sigma$
in $PGL_2(A)$.  Note that the stabilizers of the $PGL_2(A)$ action are
$$G_0 = B(A)/D(A), \qquad G_{(0,\infty)} = T(A)/D(A), \qquad G_{(0,\infty,1)} = D(A)/D(A)$$
where $B(A)$ is the upper triangular subgroup, $T(A)$ is the diagonal
subgroup, and $D(A)$ is the subgroup of scalar matrices.  The group $D(A)/D(A)$
is also the stabilizer of each orbit $(0,\infty,1,v_1,\dots ,v_{p-2})$.  
Moreover, we have $H_\bullet(B(A)/D(A),\zz) \cong H_\bullet(T(A)/D(A),\zz)$ (see \cite{nessus}).  It follows that our spectral sequence $E_{\bullet,\bullet}^\bullet(A)$ has $E^1$--term (with $\zp$ coefficients)
$$\begin{array}{|c|c|c|c|c}
                    &                     &                &       &                                    \\
H_\bullet(A^\times,\zp) & H_\bullet(A^\times,\zp) & H_\bullet(\{1\},\zp) & {\displaystyle\bop\begin{Sb}
\alpha\in A^\times \\1-\alpha \in A^\times
\end{Sb} H_\bullet(\{1\},\zp)}  &    \cdots   \\
                   &          &          &          &          \\  \hline
         \end{array}$$
and converges to $H_\bullet(PGL_2(A),\zp)$.

Now consider the cases $A=k$ and $A=R$ where
$k$ is an infinite field.  The inclusion $\pglk \ra \pgl$ induces an
injective map of spectral sequences
$$E_{\bullet,\bullet}^\bullet(k) \lra E_{\bullet,\bullet}^\bullet(R).$$

\begin{lemma}  The inclusion $\uk \ra \ur$ induces an isomorphism
$$H_\bullet(\uk,\zp) \lra H_\bullet(\ur,\zp).$$
\end{lemma}

\begin{pf}  If $G$ is an abelian group, we have an isomorphism
$$\sideset{}{_{\zp}^\bullet}\bigwedge (G\otimes \zp) \otimes \Gamma_\bullet
(\sideset{_p}{}G)
\stackrel{\cong}{\lra}  H_\bullet(G,\zp),$$
where $\sideset{_p}{}G$ is the subgroup of $G$ killed by $p$ and $\Gamma_\bullet$ is
a divided power algebra.  The lemma will follow if we can show that
$$\uk \otimes \zp \cong \ur \otimes \zp \qquad \text{and} \qquad
\sideset{_p}{}\uk \cong \sideset{_p}{}\ur.$$

We have a natural homomorphism
$$\pi: \ur \lra \uk/(\uk)^p$$
given by composing reduction modulo ${\frak m}$ with the projection
$\uk \ra \uk/(\uk)^p$.  This map is clearly surjective; we need only 
check that the kernel coincides with $(\ur)^p$.  This is an immediate
consequence of Hensel's Lemma.  Suppose $x \in \ur$ maps to 0 under $\pi$, and consider the equation $t^p-x = 0$ in $R[t]$.  Since $\pi(x) = 0$,
the reduced equation has a solution in $k$.  Since $R$ is Henselian and
$p$ is invertible in $R$, there exists a $y\in R$ with $y^p-x = 0$.  That
is, $x\in (\ur)^p$.

Note that $\sideset{_p}{}\uk$ is the subgroup of $p$th roots of unity in $k$.  A similar application of Hensel's Lemma shows that
this coincides with $\sideset{_p}{}\ur$.  This completes the proof of the lemma.
\end{pf}

\begin{cor} For each $n$, the natural map
$$H_1(\glnk,\zp) \lra H_1(\gln,\zp)$$
is an isomorphism. \hfill $\qed$
\end{cor}

\begin{cor}  The map of spectral sequences
$$E_{p,q}^1(k) \lra E_{p,q}^1(R)$$
is an isomorphism for $p\le 2$. \hfill $\qed$
\end{cor}

\begin{cor}  The natural map
$$H_2(\pglk,\zp) \lra H_2(\pgl,\zp)$$
is an isomorphism.
\end{cor}

\begin{pf}    The $E^2$--terms of the spectral
sequences look as follows:
$$E^2: \begin{array}{|cccc}
                            &     &      &             \\
        E_{0,2}^2 & E_{1,2}^2  &   0    &     \\
               0              &   0   &   0    &     \\
               \zp          &   0   &   0  &  E_{3,0}^2  \\ \hline
              \end{array}$$

By considering the commutative diagram 
$$\begin{array}{ccccccc}
 E_{3,0}^3(k) & \stackrel{d^3}{\ra} & E_{0,2}^2(k) & \ra & H_2(\pglk,\zp) &
                                      \ra &  0         \\
  \downarrow &  &  ||        &    &   \downarrow    &    &   \\
 E_{3,0}^3(R) & \stackrel{d^3}{\ra} & E_{0,2}^2(R) & 
\ra & H_2(PGL_2(R),\zp) & \ra & 0
\end{array}$$
we see that the map
$$H_2(\pglk,\zp) \lra H_2(PGL_2(R),\zp)$$
is an isomorphism.
\end{pf}

\begin{cor}  The natural map $$H_2(GL_2(k),\zp) \ra H_2(GL_2(R),\zp)$$ is an isomorphism.
\end{cor}

\begin{pf}  The Hochschild--Serre spectral sequences associated to the
extensions
$$\begin{array}{ccccccccc}
1 & \ra & \uk   & \ra & GL_2(k)  &  \ra  &  \pglk  &  \ra  &  1     \\
   &        &\downarrow&  &  \downarrow &  &  \downarrow
   &   &    \\
1 & \ra & R^\times &\ra & GL_2(R) & \ra & PGL_2(R) & \ra & 1
\end{array}$$
are isomorphic at $E_{p,q}^2$ for $0\le p \le 2$, $0 \le q \le 2$.  This
yields a commutative diagram
$$\begin{array}{cccccc}
H_3(\pglk,\zp) &\stackrel{d^3}{\ra} & E_{0,2}^3(k) & \ra & H_2(GL_2(k),\zp) & \ra\cdots \\
\downarrow &  &  ||  &  &  \downarrow  &    \\
H_3(PGL_2(R),\zp) & \stackrel{d^3}{\ra} & E_{0,2}^3(R) &
\ra & H_2(GL_2(R),\zp) & \ra\cdots 
\end{array}$$

$$\begin{array}{cccc}
\cdots\ra & H_2(\pglk,\zp) & \ra & 0 \\
      &      ||                  &      &     \\
\cdots\ra & H_2(PGL_2(R),\zp)   & \ra & 0
\end{array}$$
The result follows from the Five Lemma.
\end{pf}

\begin{cor}  The natural map $$H_2(\glnk,\zp) \ra H_2(\gln,\zp)$$
is an isomorphism for all $n$.
\end{cor}

\begin{pf}  Consider the commutative diagram
$$\begin{array}{ccc}
H_2(GL_2(k),\zp) & \stackrel{\cong}{\lra} & H_2(GL_2(R),\zp)   \\
  \downarrow  &    &  \downarrow   \\
H_2(\glnk,\zp) &  \lra  &  H_2(\gln,\zp)
\end{array}$$
The vertical arrows are isomorphisms by Suslin's stability theorem
\cite{nessus}.
\end{pf}

\noindent {\em Remark.}  The last two corollaries actually follow from
the corresponding statement for $H_2(GL)$ \cite{suslin1}.  Our approach
provides an alternate proof of this fact and also has the advantage of
proving the corresponding statement for $PGL_2$.

\medskip

\noindent {\em Remark.}  Since the map $E_{p,q}^1(k) \ra E_{p,q}^1(R)$
is an isomorphism for $p\le 2$, to show that $H_\bullet(\pglk,\zp) \cong
H_\bullet(\pgl,\zp)$ it would suffice to show that the chain complexes
$$D_\bullet(k) = E_{3,0}^1(k) \leftarrow E_{4,0}^1(k) \leftarrow \cdots$$
and
$$D_\bullet(R) = E_{3,0}^1(R) \leftarrow E_{4,0}^1(R) \leftarrow \cdots$$
are quasi-isomorphic.  One way to do this is to write down a contracting
homotopy for the quotient complex $Q_\bullet = D_\bullet(R)/D_\bullet(k)$.  Note that
the complex $Q_\bullet$ must be acyclic in the case $R=k[t_1,\dots ,t_m]/
{\frak m}^l$ for each $l\ge 2$ (see 2.1).  One would hope to be able to write down
a contracting homotopy which is independent of $l$ (or more generally,
which is independent of the fact that $R$ is a truncated polynomial ring,
but is really just a local $k$-algebra).  However, this seems to be a
very difficult question (even in the simplest case $m=1$, $l=2$).  In
fact, so far we have been unable to write down a contracting map
in the first case $Q_0 \ra Q_1$. 

\section{Bloch Groups and Third Homology}\label{ss}

We now turn our attention to the computation of $H_3(GL_2(R),\zp)$.  The
case $R=k$ was treated fully by Suslin in his beautiful paper \cite{suslin2}.  We recall the following theorem.

\begin{theorem}  There are exact sequences
$$H_3(GM_2(k),\zz) \lra H_3(GL_2(k),\zz) \lra B(k) \lra 0$$
and
$$\pi_3^0(BGM(k)^+) \lra K_3(k) \lra B(k) \lra 0.$$
\hfill $\qed$
\end{theorem}

Here, $B(k)$ denotes the Bloch group of $k$, $GM(k)$ denotes the group
of monomial matrices over $k$, and $\pi_3^0(BGM(k)^+)$ denotes the 
kernel of the canonical map
$$\pi_3(BGM(k)^+) \lra \pi_3(B\Sigma^+)$$
induced by the projection $GM(k) \ra \Sigma$ where $\Sigma$ is the
infinite symmetric group ($\Sigma = \bigcup_{n\ge 2} \Sigma_n$).

The amazing fact is that the above theorem holds for any local $k$-algebra, as we now describe.

For any local ring $A$ with infinite residue field $k$, denote by ${\cal D}(A)$ the free abelian group with basis $[x]$, where $x,1-x \in A^\times$ and define a homomorphism 
$$\phi: {\cal D}(A) \lra A^\times \otimes A^\times$$
by $\phi([x]) = x \otimes (1-x)$.  Denote by $\sigma$ the involution of
$A^\times \otimes A^\times$ given by $\sigma(x\otimes y) = -y\otimes x$, and by $(A^\times \otimes A^\times)_\sigma$ the quotient of 
$A^\times \otimes A^\times$ by the action of $\sigma$.  Define a
group ${\frak p}(A)$ by
$${\frak p}(A) = {\cal D}(A)/\langle [x] - [y] + [y/x] - [(1-x^{-1})/(1-y^{-1})] + [(1-x)/(1-y)] \rangle.$$
One checks easily (see Lemma 1.1 of \cite{suslin2}) that $\phi$ induces
a homomorphism $\overline{\phi}: {\frak p}(A) \ra (A^\times \otimes
A^\times)_\sigma$.  By definition, the {\em Bloch group}, $B(A)$, is
the kernel of $\overline{\phi}$.

Consider the action of $GL_2(A)$ on the chain complex $C_\bullet(A)$ of the
previous section.  The resulting spectral sequence is studied by Suslin in
Section 2 of \cite{suslin2}.  The main result of that section is the
existence of the exact sequence
$$H_3(GM_2(k),\zz) \lra H_3(GL_2(k),\zz) \lra B(k) \lra 0$$
(Theorem 2.1).  The proof of Theorem 2.1 goes through word for word with
$k$ replaced by $A$.  Hence we have the following result.

\begin{prop}  There is an exact sequence
$$H_3(GM_2(A),\zz) \lra H_3(GL_2(A),\zz) \lra B(A) \lra 0.$$ 
\hfill $\qed$
\end{prop}

To construct the second exact sequence, Suslin considers the action of
$GL_3(k)$ on the chain complex $C_\bullet^2(k)$, where each $C_p^2(k)$ is
the free abelian group on $(p+1)$-tuples of points in general position
in ${\Bbb P}^2(k)$.  The construction works equally well over $A$ if we
recall that closed points in ${\Bbb P}^2(A)$ consist of projective 
equivalence classes of unimodular vectors in $A^3$.  The resulting
complex $C_\bullet^2(A)$ is acyclic \cite{nessus} and the proofs of all the
results in Section 3 of \cite{suslin2} go through with one minor modification.  Lemma 3.4 of \cite{suslin2} asserts that the restriction of
the homomorphism $\partial:H_3(GL_3(k)) \ra B(k)$ to $H_3(GL_2(k))$
coincides with the homomorphism above.  The proof makes use of the
projection ${\Bbb P}^2(k) \ra {\Bbb P}^1(k)$ centered at $(0,0,1)$, which
is not well-defined for ${\Bbb P}^2(A)$.  However, a careful reading of the
proof shows that one need only construct maps among the various 
$(C_\bullet)_{GL_n(k)}$ under consideration.  The projection ${\Bbb P}^2(A)
\ra {\Bbb P}^1(A)$ {\em is} well-defined on these coinvariants since
each orbit has a unique representative for which each entry of each
vector is a unit in $A$.  This allows the proof of Lemma 3.4 to go
through for $A$.

The results in Section 4 and those through Proposition 5.1 in Section 5
may be applied word for word to $A$.  The end result is the following.

\begin{prop}[cf. Prop. 5.1, \cite{suslin2}]  There is an exact sequence
$$\pi_3^0(BGM(A)^+) \lra K_3(A) \lra B(A) \lra 0.$$
\hfill $\qed$
\end{prop}

\noindent {\em Remark.}  The rest of the proofs in Section 5 of \cite{suslin2} do not go through for $A$ as stated.  Thus, we cannot 
immediately get the exact sequence
$$0 \lra \text{Tor}(\mu(A),\mu(A))^\sim \lra K_3(A)_{\text{ind}} \lra
B(A) \lra 0.$$  However, since $\mu(R)=\mu(k)$ for a Hensel local $k$-algebra $R$ (see Lemma 3.2) and $K^M_*(R) \otimes \zp = K^M_*(k) \otimes \zp$, we see that the exact sequence exists for $R$ after 
tensoring with $\zp$ (see Proposition $4.4$ below).  Moreover, the above
results hold even more generally ({\em e.g.}, semilocal rings with infinite
residue fields).

\medskip

We now compare the exact sequences for $k$ and $R$, where $R$ is a Hensel local $k$-algebra.  Let $p$ be a prime distinct from
the characteristic of $k$.  Since tensor product is right exact, we have
a commutative diagram
$$\begin{array}{ccccccc}
\pi_3^0(BGM(k)^+)\otimes \zp& \lra& K_3(k)\otimes \zp& \lra& B(k)\otimes \zp& \lra& 0 \\
 \downarrow &  &  \downarrow &  & \downarrow &  &  \\
\pi_3^0(BGM(R)^+)\otimes \zp& \lra& K_3(R)\otimes \zp& \lra& B(R)\otimes \zp& \lra& 0.
\end{array}$$
The map $B(k) \otimes \zp \ra B(R) \otimes \zp$ is injective and
the map $K_3(k)\otimes \zp \ra K_3(R)\otimes \zp$ is an isomorphism
\cite{suslin1}.  A simple diagram chase gives the following result.

\begin{prop}  The natural map $$B(k) \otimes \zp \lra B(R) \otimes \zp$$ is an isomorphism.  \hfill $\qed$
\end{prop}

We are now in a position to compute $H_3(GL_2(R),\zp)$.

\begin{prop}  The natural map
$$H_3(GL_2(k))\otimes \zp \lra H_3(GL_2(R))\otimes \zp$$
is an  isomorphism.
\end{prop}

\begin{pf}  Consider the commutative diagram
$$\begin{array}{ccccccc}
H_3(GM_2(k))\otimes\zp & \ra & H_3(GL_2(k))\otimes\zp & \ra &
B(k) \otimes\zp & \ra & 0 \\
\downarrow &   & \downarrow &  &  \downarrow &  &  \\
H_3(GM_2(R))\otimes\zp & \ra & H_3(GL_2(R))\otimes\zp & \ra &
B(R) \otimes\zp & \ra & 0.
\end{array}$$
The map $H_3(GM_2(k))\otimes\zp \ra H_3(GM_2(R))\otimes\zp$ is
an isomorphism since we have an isomorphism
$$H_\bullet(GM_2(k),\zp) \lra H_\bullet(GM_2(R),\zp).$$
The latter isomorphism is proved by considering the Hochschild--Serre
spectral sequences associated to the extension
$$1 \lra T \lra GM_2 \lra \Sigma_2 \lra 1$$
and noting that $H_\bullet(T(k),\zp) \cong H_\bullet(T(R),\zp)$.  By the Five Lemma,
the map
$$H_3(GL_2(k))\otimes\zp \lra H_3(GL_2(R))\otimes \zp$$
is an isomorphism.
\end{pf}

\begin{cor} The natural map
$$H_3(GL_2(k),\zp) \lra H_3(GL_2(R),\zp)$$
is an isomorphism.
\end{cor}

\begin{pf} Consider the commutative diagram of universal coefficient
sequences 
$$\begin{array}{ccccc}
H_3(GL_2(k))\otimes \zp & \ra & H_3(GL_2(k),\zp) & \ra & 
\sideset{_p}{}H_2(GL_2(k))  \\
\downarrow &  &  \downarrow  &  &  \downarrow \\
H_3(GL_2(R))\otimes \zp & \ra & H_3(GL_2(R),\zp) & \ra & 
\sideset{_p}{}H_2(GL_2(R)).
\end{array}$$
The last map is an isomorphism by stability for $H_2(GL_n)$ \cite{nessus} and the corresponding statement for $H_2(GL)$.  The
result now follows from the Five Lemma.
\end{pf}

\medskip

\section{Unstable Rigidity and the Friedlander--Milnor Conjecture}\label{fmconj}

We now apply our results to the study of the Friedlander--Milnor 
conjecture in low degree.  I thank Andrei Suslin for showing me the
proofs of Propositions $5.1$ and $5.4$ below; these are well-known to
the experts.

\begin{prop}  Let $X$ be a smooth affine curve over an algebraically closed
field $k$ and let $x,y$ be closed points on $X$.  Suppose that for all $z\in X$, the inclusion $G(k)\ra G({\cal O}_z^h)$ induces an isomorphism on homology with $\zp$ coefficients.  Then the specialization homomorphisms
$$s_x,s_y: H_\bullet(G(k[X]),\zp) \lra H_\bullet(G(k),\zp)$$ coincide.
\end{prop}

\begin{pf} (Sketch)  Observe that
$$H_\bullet(G({\cal O}_x^h),\zp) = \varinjlim_{Y/X \text{\'etale}}
H_\bullet(G(k[Y]),\zp).$$
Let $\alpha \in H_\bullet(G(k[X]),\zp)$.  Denote by $\alpha(x)$ the
image of $\alpha$ under the specialization map $s_x$.  Note that since
$H_\bullet(GL_n(k),\zp)$ is a direct summand of $H_\bullet(G(k[X]),\zp)$, the
class $\alpha(x)$ also lies in the latter group.  Consider $\alpha -
\alpha(x)$.  This class has specialization $0$ and maps to $0$ in
$H_\bullet(G({\cal O}_x^h),\zp)$.  It follows that there is a curve
$Y$ which is \'etale over $X$ such that $\alpha - \alpha(x)$ maps to
$0$ in $H_\bullet(G(k[Y]),\zp)$.  Denote by $\alpha_Y$ the image of
$\alpha$ in $H_\bullet(G(k[Y]),\zp)$.  Then we have that $\alpha_Y
= \alpha(x)$, which is a constant class.  It follows that the
specialization of $\alpha_Y$ at all points of $Y$ is the same and
hence there is an affine open neighborhood of $x\in X$ on which the
specialization of $\alpha$ is constant.  Now if $y \in X$, find an
affine neighborhood of $y$ on which the specialization is constant.
Since these two neighborhoods have nontrivial intersection, we see
that $\alpha(x) = \alpha(y)$; {\em i.e.,} the homomorphisms $s_x$ and
$s_y$ agree.
\end{pf}

\begin{cor}  If $X$ is a smooth affine curve and $x,y$ are closed points on
$X$, then the specialization homomorphisms
$$s_x,s_y: H_i(GL_n(k[X]),\zp) \lra H_i(GL_n(k),\zp)$$
coincide for $i\le 3$.
\end{cor}

\begin{pf}  By Corollary $4.6$, $H_i(GL_n({\cal O}_z^h),\zp) \cong 
H_i(GL_n(k),\zp)$ for $i\le 3$.  The result follows.
\end{pf}

Recall the Friedlander--Milnor conjecture \cite{fried}.  Suppose that $G$
is a reductive group scheme over an algebraically closed field $k$.  
Denote by $BG_k$ the simplicial classifying scheme of $G$ and by $BG(k)$
the classifying space of the discrete group $G(k)$ of $k$-rational points
of $G$.

\begin{conj}  The comparison map
$$H^\bullet_{\text{\'et}}(BG_k,\zp) \lra H^\bullet(BG(k),\zp)$$
is an isomorphism.
\end{conj}

This is known to hold when $k=\overline{{\Bbb F}}_l$ \cite{fried}, and 
for the stable general linear group $GL$ \cite{suslin1}.  A few other
special cases are known \cite{sah}.

Rigidity implies the following result.

\begin{prop}  Let $k$ be an algebraically closed field and denote by
$K$ the algebraic closure of $k(T)$.  Then if rigidity holds, the natural
homomorphism
$$i_*: H_\bullet(G(k),\zp) \lra H_\bullet(G(K),\zp)$$
is an isomorphism.
\end{prop}

\begin{pf} (Sketch)  Injectivity does not require rigidity.  Indeed,
$K=\varinjlim k[C]$, where $C$ ranges over the smooth affine curves over
$k$, and each map 
$$H_\bullet(G(k),\zp)\lra H_\bullet(G(k[C]),\zp)$$
is split injective.  It follows that $i_*$ is injective (but not necessarily
split).

To prove surjectivity, note that each class in $H_\bullet(G(K),\zp)$
comes from some $\alpha \in H_\bullet(G(k[C]),\zp)$, where $C$ is a
smooth affine curve over $k$.  Denote by $C_K$ the curve
$C\times_{\text{Spec}\; k}\text{Spec}\; K$ and by $x$ the point of $C_K$
given by the inclusions $k[C] \hookrightarrow k(C) \hookrightarrow K$.
Consider the class $\alpha_K \in H_\bullet(G(K[C]),\zp)$ arising
from the homomorphism $$H_\bullet(G(k[C]),\zp) \lra
H_\bullet(G(K[C]),\zp)$$
and denote
by $\alpha_K(x)$ the image of $\alpha$ in $H_\bullet(G(K),\zp)$.
Choose a rational point $y:\text{Spec}\; k \ra C$.  This yields
another map $y:\text{Spec}\; K \ra C$.  Rigidity implies that
$\alpha_K(x) = \alpha_K(y)$ and hence $\alpha_K(x) = \alpha(y)$ (since
$y$ factors through $\text{Spec}\; k$) so
that $\alpha$ arises from a class in $H_\bullet(G(k),\zp)$. That
is, the map $i_*$ is surjective.
\end{pf}

\begin{cor}  If rigidity holds and $\text{\em char}\; k > 0$, then the
Friedlander--Milnor conjecture is true for $G$.
\end{cor}

\begin{pf}  Note that if $k$ is an algebraically closed field of characteristic $l> 0$, then $k$ is a union of transcendental extensions of
$\overline{{\Bbb F}}_l$.  It follows that the natural map
$$H_\bullet(G(\overline{{\Bbb F}}_l),\zp) \lra H_\bullet(G(k),\zp)$$
is an isomorphism.  The result follows by considering the commutative
diagram
$$\begin{array}{ccc}
H^\bullet_{\text{\'et}}(BG_k,\zp) & \lra & H^\bullet(BG(k),\zp) \\
 \downarrow  &  &  \downarrow \\
H^\bullet_{\text{\'et}}(BG_{\overline{{\Bbb F}}_l}),\zp) & \lra & H^\bullet(BG(\overline{{\Bbb F}}_l),\zp)
\end{array}$$
where the left vertical arrow is an isomorphism by the proper base change
theorem for \'etale cohomology and the bottom arrow is an isomorphism by \cite{fried}.
\end{pf}

\begin{cor}  The Friedlander--Milnor conjecture holds in positive characteristic for $H_i(GL_n)$, $i\le 3$. \hfill $\qed$
\end{cor}

{\small
\centerline{\sc Acknowledgements}

\medskip

I would like to thank Andrei Suslin for many fascinating discussions 
about rigidity and homology of congruence subgroups.  Thanks also go
to Paul Goerss and Brooke Shipley for patiently answering my questions
about cosimplicial spaces and the associated homology spectral sequence.
I also thank Philippe Elbaz-Vincent for his careful reading of the 
manuscript and for many useful comments.
Finally, I thank Mark Walker for listening to my ramblings.}

\end{document}